\documentclass[12pt]{amsart}
\usepackage{amsmath, amssymb, latexsym, amsthm}

\newcommand{\arx}[1]{}
\newcommand{\bq}{\begin{quote}}
\newcommand{\eq}{\end{quote}}
\newcommand{\cl}[1]{\overline{#1}}

\newcommand{\inv}{^{-1}}
\newcommand{\Cantor}{{{}^\N\{0,1\}}}
\newcommand{\bP}{\mathbf{P}}
\newcommand{\fU}{\mathfrak{U}}
\newcommand{\fV}{\mathfrak{V}}

\newcommand{\seq}[1]{\{#1\}_{n\in\N}}

\newcommand{\B}{\mathcal{B}}

\newcommand{\BL}{\B_\Lambda}

\newcommand{\BO}{\B_\Omega}

\newcommand{\CL}{C_\Lambda}

\newcommand{\CO}{C_\Omega}
\newcommand{\COgp}{C_{\Omega^{gp}}}
\newcommand{\CLgp}{C_{\Lambda^{gp}}}
\newcommand{\BOgp}{\B_{\Omega^{gp}}}
\newcommand{\BLgp}{\B_{\Lambda^{gp}}}

\newcommand{\cA}{\mathcal{A}}
\newcommand{\cF}{\mathcal{F}}
\newcommand{\cY}{\mathcal{Y}}

\newcommand{\N}{\naturals}
\newcommand{\NN}{{{}^{\naturals}\naturals}}

\newcommand{\cP}{\mathcal{P}}
\newcommand{\cU}{\mathcal{U}}
\newcommand{\Union}{\bigcup}
\newcommand{\cV}{\mathcal{V}}
\newcommand{\cW}{\mathcal{W}}

\newcommand{\Impl}{\Rightarrow}
\long\def\forget#1\forgotten{}

\newcommand{\rez}{\mathfrak{rez}}
\renewcommand{\b}{\mathfrak{b}}

\renewcommand{\i}{\item}

\renewcommand{\u}{\mathfrak{u}}

\newcommand{\w}{\omega}

\newcommand{\nin}{\not\in}

\newcommand{\sbst}{\subseteq}

\newcommand{\sm}{\setminus}

\renewcommand{\pi}{pseudo-intersection}

\newcommand{\non}{\mathsf{non}}

\newtheorem{thm}{Theorem}

\newtheorem{prob}[thm]{Problem}

\newtheorem{cor}[thm]{Corollary}
\newtheorem{conj}[thm]{Conjecture}
\theoremstyle{definition}

\theoremstyle{remark}

\newcommand{\be}{\begin{enumerate}}
\newcommand{\ee}{\end{enumerate}}
\newcommand{\bi}{\begin{itemize}}
\newcommand{\ei}{\end{itemize}}

\newcommand{\sfin}{\mathsf{S}_{fin}}

\renewcommand{\split}{\mathsf{Split}}

\newcommand{\naturals}{{\mathbb N}}

\author{Boaz Tsaban}
\thanks{Partially supported by the Golda Meir Fund and
the Edmund Landau Center for Research in Mathematical Analysis and Related Areas,
sponsored by the Minerva Foundation (Germany).}
\address{Einstein Institute of Mathematics, Hebrew University of Jerusalem,
Givat Ram, Jerusalem 91904, Israel}
\email{tsaban@math.huji.ac.il}
\urladdr{http://www.cs.biu.ac.il/\~{}tsaban}

\title[The Reznichenko property]{The minimal cardinality where the Reznichenko property fails}

\begin{document}
\begin{abstract}
A topological space $X$ has the Fr\'echet-Urysohn
property if for each subset $A$ of $X$ and each element
$x$ in $\cl{A}$, there exists a countable sequence
of elements of $A$ which converges to $x$.
Reznichenko introduced a natural generalization of
this property,
where the converging sequence of elements is replaced by a
sequence of disjoint finite sets which eventually intersect
each neighborhood of $x$.
In \cite{FunRez}, Ko\v{c}inac and Scheepers conjecture:
\begin{quote}
The minimal cardinality of a set $X$ of real numbers
such that $C_p(X)$ does not have the weak Fr\'echet-Urysohn property
is equal to $\b$.
\end{quote}
($\b$ is the minimal cardinality of an unbounded family in the
Baire space $\NN$).
We prove the Ko\v{c}inac-Scheepers conjecture by showing that
if $C_p(X)$ has the Reznichenko property, then a continuous image of
$X$ cannot be a subbase for a non-feeble filter on $\N$.
\end{abstract}


\maketitle

\section{Introduction}

A topological space $X$ has the Fr\'echet-Urysohn
property if for each subset $A$ of $X$ and each
$x\in\cl{A}$, there exists a sequence $\seq{a_n}$
of elements of $A$ which converges to $x$.
If $x\nin A$ then we may assume that the elements
$a_n$, $n\in\N$, are distinct.
The following natural generalization of this property
was introduced by Reznichenko \cite{MaTi}:
\bq
For each subset $A$ of $X$ and each element
$x$ in $\cl{A}\sm A$, there exists a countably infinite
pairwise disjoint collection $\cF$ of finite subsets of $A$
such that for each neighborhood $U$ of $x$, $U\cap F\neq\emptyset$ for
all but finitely many $F\in\cF$.
\eq
In \cite{MaTi} this is called the \emph{weak Fr\'echet-Urysohn} property.
In other works \cite{FunRez, coc7, SakaiRez}
this also appears as the \emph{Reznichenko} proeprty.

For a topological space $X$ denote by $C_p(X)$ the space of continuous
real-valued functions with the topology of pointwise convergence.
A comprehensive duality theory was developed by Arkhangel'ski\v{i}
and others (see, e.g., \cite{GN, Sakai, FunRez, coc7} and references therein)
which characterizes topological properties of $C_p(X)$ for a Tychonoff space
$X$ in terms of covering properties of $X$.
In \cite{FunRez, coc7} this is done for a conjunction of the
Reznichenko property and some other classical property
(countable strong fan tightness in \cite{FunRez} and
countable fan tightness in \cite{coc7}).
According to Sakai \cite{Sakai}, a space $X$ has countable fan tightness
if for each $x\in X$ and each sequence $\seq{A_n}$ of subsets of $X$
with $x\in\cl{A_n}\sm A_n$ for each $n$, there exist finite sets $F_n\sbst A_n$,
$n\in\N$, such that $x\in \cl{\Union_nF_n}$.
In Theorem 19 of \cite{coc7}, Ko\v{c}inac and Scheepers prove that
for a Tychonoff space $X$,
$C_p(X)$ has countable fan tightness as well as Reznichenko's property if, and only if,
each finite power of $X$ has the Hurewicz covering property.

The \emph{Baire space} $\NN$ of infinite sequences
of natural numbers is equipped with the product topology
(where the topology of $\N$ is discrete).
A quasiordering $\le^*$ is defined on the Baire space $\NN$ by eventual
dominance:
$$f\le^* g\quad\mbox{if}\qquad f(n)\le g(n)\mbox{ for all but finitely many }n.$$
We say that a subset $Y$ of $\NN$ is \emph{bounded} if
there exists $g$ in $\NN$ such that for each $f\in Y$, $f\le^* g$.
Otherwise, we say that $Y$ is \emph{unbounded}.
$\b$ denotes the minimal cardinality of an unbounded family in $\NN$.
According to a theorem of Hurewicz
\cite{HURE27},
a set of reals $X$ has the Hurewicz property if, and only if,
each continuous image of $X$ in $\NN$ is bounded.
This and the preceding discussion imply that for each set of reals $X$
of cardinality smaller than $\b$, $C_p(X)$ has the Reznichenko property.
Ko\v{c}inac and Scheepers conclude their paper \cite{FunRez}
with the following.
\begin{conj}
$\b$ is the minimal cardinality of a set $X$ of real numbers
such that $C_p(X)$ does not have the Reznichenko property.
\end{conj}
We prove that this conjecture is true.

\section{A proof of the Ko\v{c}inac-Scheepers conjecture}

Throughout the paper, when we say that
$\cU$ is a \emph{cover} of $X$ we mean that $X\sbst\cup\cU$ but
$X$ is not contained in any member of $\cU$.
A cover $\cU$ of a space $X$ is an \emph{$\w$-cover} of $X$
if each finite subset $F$ of $X$ is contained in some member of
$\cU$. This notion is due to Gerlits and Nagy \cite{GN},
and is starring in \cite{FunRez, coc7}.
According to \cite{FunRez, coc7}, a cover $\cU$ of $X$
is \emph{$\w$-groupable} if there exists a partition $\cP$ of
$\cU$ into finite sets such that for each finite $F\sbst X$ and
all but finitely many $\cF\in\cP$,
there exists $U\in\cF$ such that $F\sbst U$.
Thus, each \emph{$\w$-groupable} cover is an $\w$-cover and
contains a countable $\w$-groupable cover.

In \cite{coc7} it is proved that if each open $\w$-cover of a set of reals
$X$ is $\w$-groupable and $C_p(X)$ has countable fan tightness, then
$C_p(X)$ has the Reznichenko property.
Recently, Sakai \cite{SakaiRez} proved that the assumption of
countable fan tightness is not needed here. More precisely,
say that an open $\w$-cover $\cU$ of $X$ is
\emph{$\w$-shrinkable} if for each $U \in \cU$ there exists
a closed subset $C_U\sbst U$ such that
$\{C_U : U \in \cU\}$ is an $\w$-cover of $X$.
Then the following duality result holds.
\begin{thm}[Sakai \cite{SakaiRez}]\label{sakaithm}
For a Tychonoff space $X$, the following are equivalent:
\be
\i $C_p(X)$ has the Reznichenko property;
\i Each $\w$-shrinkable open $\w$-cover of $X$
is $\w$-groupable.
\ee
\end{thm}
It is the other direction of this result that we are interested in here.
Observe that any \emph{clopen} $\w$-cover is trivially $\w$-shrinkable.
\begin{cor}\label{w-chs-gp}
Assume that $X$ is a Tychonoff space such that $C_p(X)$ has the Reznichenko property.
Then each clopen $\w$-cover of $X$ is $\w$-groupable.
\end{cor}

From now on $X$ will always denote a set of reals.
As all powers of sets of reals are Lindel\"of,
we may assume that all covers we consider are
countable \cite{GN}.
For conciseness, we introduce some notation.
For collections of covers of $X$ $\fU$ and $\fV$,
we say that $X$ satisfies $\binom{\fU}{\fV}$ (read: $\fU$ choose $\fV$)
if each element of $\fU$ contains an element of $\fV$ \cite{tautau}.
Let $\CO$ and $\COgp$ denote the collections of clopen
$\w$-covers and clopen $\w$-groupable covers of $X$,
respectively.
Corollary \ref{w-chs-gp} says that the Reznichenko property for $C_p(X)$
implies $\binom{\CO}{\COgp}$.

As a warm up towards the real solution, we make the following observation.
According to \cite{coc1}, a space $X$ satisfies $\split(\fU,\fV)$ if
every cover $\cU\in\fU$ can be split into two disjoint subcovers $\cV$ and $\cW$
which contain elements of $\fV$.
Observe that $\binom{\CO}{\COgp}$ implies $\split(\CO,\CO)$.
The \emph{critical cardinality} of a property $\bP$ (or collection) of sets of reals, $\non(\bP)$,
is the minimal cardinality of a set of reals which does not satisfy this property.
Write
$$\rez=\non(\{X : C_p(X)\mbox{ has the Reznichenko property}\}).$$
Then we know that $\b\le\rez$, and
the Ko\v{c}inac-Scheepers conjecture asserts that
$\rez=\b$. By Corollary \ref{w-chs-gp},
$\rez\le\non(\split(\CO,\CO))$.
In \cite{coc2} it is proved that $\non(\split(\CO,\CO))=\u$,
where $\u$ is the \emph{ultrafilter number}
denoting the minimal size of a base for a nonprincipal ultrafilter on $\N$.
Consequently, $\rez\le\u$. It is well known that $\b\le\u$, but it is
consistent that $\b<\u$. Thus this does not prove the conjecture.
However, this is the approach that we will use:
We will use the language of filters to prove
that $\non(\binom{\CO}{\COgp})=\b$. By Corollary \ref{w-chs-gp},
$\b\le\rez\le\non(\binom{\CO}{\COgp})$, so this will suffice.

A \emph{nonprincipal filter} on $\N$ is a family $\cF\sbst P(\N)$
that contains all cofinite sets but not the empty set,
is closed under supersets, and
is closed under finite intersections
(in particular, all elements of a nonprincipal filter are infinite).
A \emph{base} $\B$ for a nonprincipal filter $\cF$ is a subfamily of $\cF$
such that for each $A\in\cF$ there exists $B\in\B$ such that
$B\sbst A$. If the closure of $\B$ under finite intersections is
a base for a nonprincipal filter $\cF$, then we say that $\B$ is a
\emph{subbase} for $\cF$.
A family $\cY\sbst P(\N)$ is \emph{centered} if for each
finite subset $\cA$ of $\cY$, $\cap\cA$ is infinite.
Thus a subbase $\B$ for a nonprincipal filter is a centered
family such that
for each $n$ there exists $B\in\B$ with $n\nin B$.
For a nonprincipal filter $\cF$ on $\N$ and a finite-to-one function
$f:\N\to\N$, $f(\cF):=\{A\sbst\N : f\inv[A]\in\cF\}$ is again a nonprincipal
filter on $\N$.

A filter $\cF$ is \emph{feeble} if
there exists a finite-to-one function $f$ such that
$f(\cF)$ consists of only the cofinite sets.
$\cF$ is feeble if, and only if, there exists a partition
$\seq{F_n}$ of $\N$ into finite sets such that for each
$A\in\cF$, $A\cap F_n\neq\emptyset$ for all but finitely many $n$
(take $F_n = f\inv[\{n\}]$). Thus $\B$ is a subbase for a feeble filter
if, and only if:
\be
\i $\B$ is centered,
\i For each $n$ there exists $B\in\B$ such that $n\nin B$; and
\i There exists a partition
$\seq{F_n}$ of $\N$ into finite sets such that
for each $k$ and each $A_1,\dots,A_k\in\B$,
$A_1\cap\dots\cap A_k\cap F_n\neq\emptyset$
for all but finitely many $n$.
\ee

Define a topology on $P(\N)$ by identifying it with \emph{Cantor's space}
$\Cantor$ (which is equipped with the product topology).
\begin{thm}\label{feeble}
For a set of reals $X$, the following are equivalent:
\be
\i $X$ satisfies $\binom{\CO}{\COgp}$;
\i For each continuous function $\Psi:X\to P(\N)$,
$\Psi[X]$ is not a subbase for a non-feeble filter on $\N$.
\ee
\end{thm}
\begin{proof}
$(1\Impl 2)$ Assume that $\Psi:X\to P(\N)$ is continuous
and $\B=\Psi[X]$ is a subbase for a nonprincipal filter $\cF$
on $\N$.
Consider the (clopen!) subsets $O_n = \{A\sbst\N : n\in A\}$, $n\in\N$,
of $P(\N)$. For each $n$, there exists $B\in\B$
such that $B\nin O_n$ ($n\nin B$), thus $X\not\sbst\Psi\inv[O_n]$.

As $\B$ is centered, $\seq{O_n}$ is an $\w$-cover of $\B$,
and therefore $\seq{\Psi\inv[O_n]}$ is a clopen $\w$-cover of $X$.
Let $A\sbst\N$ be such that the enumeration $\{\Psi\inv[O_n]\}_{n\in A}$
is bijective.
Apply $\binom{\CO}{\COgp}$ to obtain a partition $\seq{F_n}$ of
$A$ into finite sets such that for each finite $F\sbst X$,
and all but finitely many $n$, there exists $m\in F_n$ such that
$F\sbst\Psi\inv[O_m]$ (that is, $\Psi[F]\sbst O_m$, or
$\bigcap_{x\in F}\Psi(x)\cap F_n\neq\emptyset$).
Add to each $F_n$ an element from $\N\sm A$ so that $\seq{F_n}$
becomes a partition of $\N$. Then the sequence $\seq{F_n}$
witnesses that $\B$ is a subbase for a \emph{feeble} filter.

$(2\Impl 1)$ Assume that $\cU=\seq{U_n}$ is a clopen $\w$-cover
of $X$. Define $\Psi:X\to P(\N)$ by
$$\Psi(x) = \{n : x\in U_n\}.$$
As $\cU$ is clopen, $\Psi$ is continuous.
As $\cU$ is an $\w$-cover of $X$,
$\B=\Psi[X]$ is centered (see Lemma 2.2 in \cite{tau}).
For each $n$ there exists $x\in X\sm U_n$, thus
$n\nin\Psi(x)$.
Therefore $\B$ is a subbase for a feeble filter. Fix
a partition $\seq{F_n}$ of $\N$ into finite sets such that
for each $\Psi(x_1),\dots,\Psi(x_k)\in\B$,
$\Psi(x_1)\cap\dots\cap\Psi(x_k)\cap F_n\neq\emptyset$
(that is, there exists $m\in F_n$ such that $x_1,\dots,x_k\in U_m$)
for all but finitely many $n$.
This shows that $\cU$ is groupable.
\end{proof}

\begin{cor}
$\non(\binom{\CO}{\COgp})=\b$.
\end{cor}
\begin{proof}
Every nonprincipal filter on $\N$ with a
(sub)base of cardinality smaller than $\b$ is feeble
(essentially, \cite{Solomon}), and by an unpublished result of
Petr Simon, there exists a non-feeble filter with a (sub)base of cardinality $\b$
-- see \cite{BlassHBK} for the proofs. Use Theorem \ref{feeble}.
\end{proof}

This completes the proof of the Ko\v{c}inac-Scheepers conjecture.

\section{Consequences and open problems}
Let $\BO$ and $\BOgp$ denote the collections of \emph{countable Borel}
$\w$-covers and $\w$-groupable covers of $X$, respectively.
The same proof as in Theorem \ref{feeble} shows that the analogue
theorem where ``continuous'' is replaced by ``Borel'' holds.

$\cU$ is a \emph{large} cover of a space $X$ if each member of $X$ is
contained in infinitely many members of $\cU$.
Let $\BL$, $\Lambda$, and $\CL$ denote the collections of countable large
Borel, open, and clopen covers of $X$, respectively.
According to \cite{coc7}, a large cover $\cU$ of $X$ is \emph{groupable}
if there exists a partition $\cP$ of $\cU$ into finite sets such that
for each $x\in X$ and all but finitely many $\cF\in\cP$, $x\in\cup\cF$.
Let $\BLgp$, $\Lambda^{gp}$, and $\CLgp$ denote the collections of countable groupable
Borel, open, and clopen covers of $X$, respectively.

\begin{cor}
The critical cardinalities of the classes
$\binom{\BL}{\BLgp}$,
$\binom{\BO}{\BOgp}$,
$\binom{\BO}{\BLgp}$,
$\binom{\Lambda}{\Lambda^{gp}}$,
$\binom{\Omega}{\Omega^{gp}}$,
$\binom{\Omega}{\Lambda^{gp}}$,
$\binom{\CL}{\CLgp}$,
$\binom{\CO}{\COgp}$, and
$\binom{\CO}{\CLgp}$
are all equal to $\b$.
\end{cor}
\begin{proof}
By the Borel version of Theorem \ref{feeble},
$\non(\binom{\BO}{\BOgp})=\b$.
In \cite{hureslaloms} it is proved that $\non(\binom{\BL}{\BLgp})=\b$.
These two properties imply all other properties in the list.
Now, all properties in the list imply either $\binom{\CL}{\CLgp}$
or $\binom{\CO}{\CLgp}$, whose critical cardinality is $\b$ by
Theorem \ref{feeble} and \cite{hureslaloms}.
\end{proof}

If we forget about the topology and consider \emph{arbitrary} countable covers,
we get the following characterization of $\b$, which extends
Theorem 15 of \cite{coc7} and Corollary 2.7 of \cite{hureslaloms}.
For a cardinal $\kappa$, denote by
$\Lambda_\kappa$, $\Omega_\kappa$, $\Lambda^{gp}_\kappa$, and $\Omega^{gp}_\kappa$
the collections of countable large covers, $\w$-covers, groupable covers, and $\w$-groupable
covers of $\kappa$, respectively.
\begin{cor}\label{bchar}
For an infinite cardinal $\kappa$, the following are equivalent:
\be
\i $\kappa<\b$,
\i $\binom{\Lambda_\kappa}{\Lambda^{gp}_\kappa}$,
\i $\binom{\Omega_\kappa}{\Lambda^{gp}_\kappa}$; and
\i $\binom{\Omega_\kappa}{\Omega^{gp}_\kappa}$.
\ee
\end{cor}

It is an open problem \cite{SakaiRez} whether
item (2) in Sakai's Theorem \ref{sakaithm}
can be replaced with $\binom{\Omega}{\Omega^{gp}}$
(by the theorem, if $X$ satisfies $\binom{\Omega}{\Omega^{gp}}$,
then $C_p(X)$ has the Reznichenko property; the other direction is the unclear
one).

For collections $\fU$ and $\fV$ of covers of $X$, we say that
$X$ satisfies $\sfin(\fU,\fV)$ if:
\bq
For each sequence $\seq{\cU_n}$
of members of $\fU$, there is a sequence
$\seq{\cF_n}$ such that each $\cF_n$ is a finite
subset of $\cU_n$, and $\Union_{n\in\N}\cF_n\in\fV$.
\eq
In \cite{hureslaloms} it is proved that
$\binom{\Lambda}{\Lambda^{gp}}=\sfin(\Lambda,\Lambda^{gp})$
(which is the same as the Hurewicz covering property \cite{coc7}).
We do not know whether the analogue result for $\binom{\Omega}{\Omega^{gp}}$ is true.
\begin{prob}\label{w/wgp}
Does $\binom{\Omega}{\Omega^{gp}}=\sfin(\Omega,\Omega^{gp})$?
\end{prob}

In \cite{coc7} it is proved that $X$ satisfies $\sfin(\Omega,\Omega^{gp})$
if, and only if, all finite powers of $X$ satisfy the Hurewicz
covering property $\sfin(\Lambda,\Lambda^{gp})$, which we now know is the
same as $\binom{\Lambda}{\Lambda^{gp}}$.

\bigskip

\emph{Added after publication.}
\textsl{The answer to Problem \ref{w/wgp} is
``No'', in the following strong sense: Masami Sakai
proves in: 
\emph{Two properties of $C_p(X)$ weaker than the Fr\'echet Urysohn property},
Topology and its Applications \textbf{153} (2006), 2795--2804,
that every analytic set of reals
(and, in particular, the Baire space $\NN$) satisfies $\binom{\BO}{\BOgp}$.
But we know that $\NN$ does not satisfy the Hurewicz covering property.
}

\end{document}